\documentstyle[11pt]{article}
\newtheorem{definition}{Definition}
\newtheorem{theorem}{Theorem}
\newtheorem{remark}{Remark}
\newtheorem{lemma}{Lemma}

\newtheorem{corollary}{Corollary}

\newcommand{\bC}{{\bf C}}

\newcommand{\bR}{{\bf R}}

\newcommand{\la}{\lambda}

\newcommand{\si}{\sigma}

\newcommand{\comb}[2]{\pmatrix{#1\cr #2}}
\newcommand{\om}{\omega}
\newcommand{\Ga}{\Gamma}
\newcommand{\Lam}{\Lambda}

\newcommand{\Spin}{\mbox{\rm Spin}}

\newcommand{\pul}{{{\frac{1}{2}}}}

\newcommand{\reo}{\mbox{${\bf R_0^{n+1}}$}}

\newcommand{\ps}{\mbox{${\bf S}_{\frac{1}{2}}$}}

\begin{document}
\title{The higher spin Dirac operators.\thanks{The research was
supported by grant GAUK 713. This paper is in final form and no version
of it will be published elsewhere.}}
\author{Jarol\'{\i}m Bure\v{s}, Prague}
\date{}
\maketitle
\begin{abstract}
There is a certain family of conformally invariant first order
elliptic systems which include the Dirac operator as its first
and simplest member.
Their general definition is given and some of their basic properties
are described.
 A special attention
is paid to the Rarita-Schwinger operator, the second simplest operator
in the row. Its basic properties are
described in more details. In the last
part  indices of discussed operators are computed.


\noindent MSC 1991 : 58 G 25, 53 C 25.
\end{abstract}

\section{Introduction}

Let $M$ be a smooth oriented compact n-dimensional
manifold,
endowed with a Riemannian metric and a spin structure.
A lot of information has been collected concerning
basic invariant differential operators
on $M$ as the Laplace (of the second order) and the Dirac operator
(of the first order).
These invariant operators and their properties are strongly
related with geometry and topology of $M$.

Their behavior under conformal change of metric allowed us to
study their conformal invariance and put it into general scheme
of invariant operators on AHS manifolds.
There is a list of conformal invariant operators of different
order, all of them are defined on any Riemannian spin manifold.
For this operators similar notions can be defined and properties and relations
as for
basic operators (Dirac and Laplace) can be studied.

Recently, a growing interest is paid to properties of more
complicated  invariant first order differential operators on
$M.$ A prototype of them is the Rarita-Schwinger
operator (see
\cite{Fra1,Fra2,FraS,MP,Pe1,Pe2,Pe3,Pe4,Pe5,RaS,US,Wa}).
It acts on sections of the bundle associated
to a more complicated representation of the group $\Spin(n).$

It it also included in a series of operators on $M$ which are
called higher spin Dirac operators  and also in a series of
operators called higher Rarita-Schwinger operators which are
defined on spin valued symmetric tensor fields
and has been intensively used also in
Clifford analysis in connection with monogenic differential
forms (see \cite{DSS,So1,So2,SoS,So3}) or symmetric functions
(see\cite{Som1}).

In general, there  is a question
which new geometrical and topological
characterizations of
Riemannian spin manifold can be obtained from properties of this
additional
invariant operators on it.
To try to see an answer, it is first necessary to learn more
about properties of these operators.
The aim of the paper is to collect and review facts which
exist  mostly only in a preprint form and
to add some new facts (index computations).
In the contribution,
the operators are defined including their normalisation
and some of their properties are described.
 Main attention is
concentrated to the Rarita-Schwinger operator.
Basic questions discussed in the paper are:

1. The description of conformally invariant first order differential operators
(including the Rarita-Schwinger operator).
\vskip 2mm
2. The spectrum of Rarita Schwinger operator on the flat model, i.e. on
spheres.
\vskip 2mm
3. A complete description of polynomial solutions of the
Rarita-Schwinger equation.
\vskip 2mm
4. The index of Rarita Schwinger operator and higher spin Dirac
operators.

The first three parts have a review character, most results there
are taken from papers which are at present only
in a preprint form, the index properties are new results.

\section{First order conformally invariant operators.}

There is a scheme for a construction of conformal invariant operators
let us recall it shortly from (\cite{BuSo}).

Let M be a compact oriented spin manifold with a conformal structure.
Fix a  Riemannian metric $g$ in the given conformal class
then we have on $M$ principal
fibre bundles
$$
\tilde{\cal P}\equiv\tilde{\cal P}_{Spin}\rightarrow
{\cal P}_{SO} \rightarrow M.
$$

Finite-dimensional irreducible representations ${\bf V}_{\la}$ of the group
$\Spin(n)$ are determined by their highest weights $\la \in
\Lambda^{+},$
                    where for $n=2k$ even, we have
$$
\Lambda^+ =
\{ \la =(\la_1,...,\la_k); \la_1\geq \la_2\geq ... \geq \la_{k-1}\geq
|\la_k|\},
\la_i \in {\bf Z}\cup \pul {\bf Z}
$$
 and for $n=2k+1$ odd, we have
$$
\Lambda^+ =
\{ \la =(\la_1,...,\la_k); \la_1\geq \la_2\geq ... \geq \la_{k-1}\geq
\la_k\geq 0\},
\la_i \in {\bf Z}\cup \pul {\bf Z}.
$$

Invariant operators are acting among spaces of sections of
the corresponding associated bundles
$$
V_{\la} = {\tilde{\cal P}}\times_{\Spin(n)} {\bf V}_{\la}
$$
over $M.$
%\vskip 1mm
Let us consider the Levi-Civita connection $\omega$ of the chosen
Riemannian metric on ${\cal P}$ and let $\tilde{\omega}$ be its (unique) lift
to ${\tilde{\cal P}}.$
For any choice of ${\la}\in\Lam^+,$ we have the associated covariant derivative
$$
\nabla_{\la}:\Gamma(V_{\la})
{\rightarrow} \Gamma(V_{\la}
\otimes T^*(M)).
$$
Tensor
product ${\bf V_{\la}}\otimes {\bf C}_n$ can be decomposed
into irreducible components
$$
{\bf V_{\la}}\otimes {\bf C}_n = \oplus_{\la'\in A} {\bf V}_{\la'},
$$
where $A$ is the set of  highest weights of all irreducible
components (multiplicities included).
There are simple rules how to describe $A=A(\la)$ explicitely for any $\la$
(see \cite{F}).
Let $\pi_{\la'}$ be the projection from
${\bf V}_{\la}\otimes {\bf C}_n$ to
${\bf V}_{\la'}.$ Then  operators
$$
D_{\la,\la'} :\Gamma(V_{\la}) \rightarrow \Gamma (V_{\la'}),\;
D_{\la,\la'}:=\pi_{\la'}\circ\nabla^{\la}
$$
are first order conformally invariant differential operators
and all such operators
can be  constructed in this  way.

Any conformally invariant first order differential
 operator is uniquely determined
(up to a constant multiple) by a choice of allowed $\la$ and $\la'$
but there is no natural normalization in general.
 To study spectral properties, it is necessary fix a scale of the operator, to choose appropriate
 normalization. For the Dirac operator, the choice of normalization
 is given by the Clifford action. By using twisted Dirac operators,
 we shall extend this normalization to a wide class of first order
 operators (which includes our higher spin Dirac perators as well as
higher Rarita-Schwinger operators) .

\begin{definition}\cite{BuSo}
 Let ${\bf S}_{\pul}^{}$ (for $n=2k+1$), resp.
 ${\bf S}_{\pul}^{}={\bf S}_{\pul}^{+}\oplus{\bf S}_{\pul}^{-}$
(for $n=2k$), denote the basic spinor representations
 with highest weights
$\si=({\frac{1}{2}},\ldots,{\frac{1}{2}},{\frac{1}{2}}),$
resp.
$\si^{\pm}=({\frac{1}{2}},\ldots,{\frac{1}{2}},\pm{\frac{1}{2}}).$

Let $\la
\in\Lam^+,$ (for $n=2k+1$), resp. $\la^{\pm}\in\Lam^+$ (for
$n=2k$)
 be  dominant weights with
$\la =(\la_1,...,\la_{k-1},{\frac{1}{2}}),$
resp. $\la^{\pm}=(\la_1,...,\la_{k-1},\pm{\frac{1}{2}}).$
 Denote further $\la'=\la-\si\in\Lam^+,$ resp.
${\la'}=\la^{+}-\si^{+}\in\Lam^+.$
In even dimensions, we shall use the notation
$$
{\bf V}_{\la}={\bf V}_{\la^+}\oplus {\bf V}_{\la^-}.
$$

The representation ${\bf V}_{\la}$  appears with multiplicity one
in the decomposition of the tensor product
${\bf S}_{\pul}\otimes {\bf V}_{\la'}$ (it is the Cartan product of both
representations). Hence we can write the product
as
$$
{\bf S}_{\pul}\otimes {\bf V}_{\la'}=
{\bf V}_{\la} \oplus {\bf W},
$$
where ${\bf W}$ is the sum of all other  irreducible components
in the decomposition.

Let  $D^T_{\la'}$ be the twisted Dirac operator
on $S_{\pul}\otimes V_{\la'}.$  If we write the operator $D^T_{\la'}$
in the block form as

\vskip 1.7cm
{\center{
%TexCad Options
%\grade{\off}
%\emlines{\off}
%\beziermacro{\on}
%\reduce{\on}
%\snapping{\off}
%\quality{2.00}
%\graddiff{0.01}
%\snapasp{1}
%\zoom{1.50}
\unitlength 1.00mm
\linethickness{0.4pt}
\begin{picture}(73,50)

\put(25,57){$\Gamma(S_{\pul}\otimes V_{\lambda'})$}
\put(60,57){$\Gamma(S_{\pul}\otimes V_{\lambda'})$}

\put(30,42){$\Gamma(V_{\lambda})$}
\put(65,42){$\Gamma(V_{\lambda})$}

\put(30,27){$\Gamma(W)$}
\put(65,27){$\Gamma(W)$}

\put(46,58){\vector(1,0){11.5}}
\put(42,43){\vector(1,0){17.5}}
\put(42,28){\vector(1,0){17.5}}
\put(42,31.3){\vector(2,1){17}}
\put(42,39.7){\vector(2,-1){17}}

\put(32,50){$\|$}
\put(67,50){$\|$}

\put(32,35){$\oplus$}
\put(67,35){$\oplus$}

\put(49,62){\makebox(0,0)[cc]{$D^T_{\lambda'}$}}
\put(49,47){\makebox(0,0)[cc]{$D_{\lambda}$}}
%\put(48,30){\makebox(0,0)[cc]{$D_{\lambda}^{\perp}$}}

\put(16,52){\makebox(0,0)[cc]{$$}}
\put(40,52){\makebox(0,0)[cc]{$$}}

\end{picture}
}}

\vskip -1.7cm
\noindent
we have defined four invariant operators, one of them being the operator
$$
D_{\la}:\Ga(V_{\la})\rightarrow\Ga(V_{\la}).
$$
Operators $D_{\la}$ defined in such a way will be called
{\it generalized (higher spin) Dirac operators}.
\end{definition}
\vskip 2mm
The case of Rarita-Schwinger operator is included at the beginning of the
scheme, as follows.

\vskip 2mm

A certain subclass of invariant operators
discussed above is related with the following
higher dimensional generalizations of holomorphic differential
forms (see \cite{DSS,So2}).
Spinor valued differential forms are
elements of the twisted de Rham sequence (\cite{BuSo}):

$$  \Gamma(S_{\pul}^{\pm})\stackrel{\nabla^S}{\rightarrow}
%  \Gamma({\Omega}^{1}_c\otimes S_{\pul}^{\pm})\stackrel{\nabla^S}{\rightarrow} &
  \ldots
  \Gamma({\Omega}^{k}_c \otimes S_{\pul}^{\pm})\stackrel{\nabla^S}{\rightarrow}
  \ldots
 % \Gamma({\Omega}^{2k-1}_c \otimes S_{\pul}^{\pm})
\stackrel{\nabla^S}{\rightarrow}
  \Gamma({\Omega}^{n}_c \otimes S_{\pul}^{\pm})
$$
where $\nabla^S$ denotes the associated covariant  derivative
on spinor bundles extended to $S_{\pul}$-valued forms (see
\cite{So2,VSe}).

Every representation
$\Lambda^k({\bf C}_n)\otimes {\bf S}_{\pul}^{}$  can be split
into irreducible pieces. There are no multiplicities in the
decomposition, so the irreducible pieces are well defined.
For $k$ forms ($k\leq[n/2]$), there are  $k$ pieces in the
decomposition and the decomposition is symmetric with
respect to the action of the Hodge star operator.
The space of spinor valued $k$-forms
$ \Gamma({\Omega}^{k}_c \otimes S^{\pm})$
($k\leq [n/2]$) can be written
as the sum $\oplus_{j=1}^kE^{k,j}$
and it can be checked (see \cite{DSS,VSe,So2}) that
$E^{k,j}$ is the bundle associated with the representation
with the highest weight
$\la_j=({\frac{3}{2}},\ldots,{\frac{3}{2}},{\frac{1}{2}},\ldots,
{\frac{1}{2}},
\pm{\frac{1}{2}}),$ where the number $j$ indicates that
the component ${\frac{3}{2}}$ appears with multiplicity equal to $j.$
Signs $\pm$ at the last components are relevant only  in even
dimensions (more details can be found in \cite{VSe}).
The whole splitting can be described by the following triangle
shaped diagram (in odd dimensions, there are two columns of the
same length in the middle).

%\hskip -7mm
\begin{displaymath}
\begin{array}{ccccccccccccc}
  E^{0,0} &\stackrel{D_0}{\longrightarrow}&
  E^{1,0} &\stackrel{D_0}{\longrightarrow}&
  \ldots &     \stackrel{D_0}{\longrightarrow}&
  E^{k,0} &\stackrel{D_0}{\longrightarrow}&
  \ldots & \stackrel{D_0}{\longrightarrow}&
  E^{2k-1,0} &\stackrel{D_0}{\longrightarrow}&
  E^{2k,0} \\

  && \oplus && \oplus && \oplus && \oplus && \oplus &&  \\

  &&
  E^{1,1} &\stackrel{D_1}{\longrightarrow}&
  \ldots &          \stackrel{D_1}{\longrightarrow}&
  E^{k,1} &     \stackrel{D_1}{\longrightarrow}&
  \ldots & \stackrel{D_1}{\longrightarrow}&
  E^{2k-1,1}
  &&
   \\
   &&  &&   \oplus && \oplus && \oplus &&  && \\
   &&  && \ldots &\stackrel{D_{j}}{\longrightarrow}& \vdots
   &\stackrel{D_{j}}{\longrightarrow}& \ldots &&  && \\
    &&  &&  && \oplus &&  &&  && \\
    &&  &&  && E^{k,k} &&  &&  && \\
\end{array}
\end{displaymath}

The general construction of invariant operators described above
can be used
in the special case of spinor valued forms.
The  covariant
derivative $\nabla ^S$ restricted to $E^{k,j}$ and projected
to $E^{k+1,j'}$ is an example of this general construction.
It can be shown that if $|j-j'|>1,$ then the corresponding
invariant operator is trivial. We shall be mainly interested
in 'horizontal arrows', i.e. in operators $D_j$ given
by restriction to $E^{k,j}$ and projection to $E^{k+1,j}.$
They are indicated in the above scheme.

The simplest cases among
them are well known. The operator $D_0$ is (a multiple of) the Dirac
operator.

The operator $D_1$ is (an elliptic version of)
the operator called Rarita-Schwinger operator by physicists
(see \cite{RaS,Wa}).

All of them are elliptic operators  (see
\cite{So1, Br}).

Note that all operators $D_j$ on the same row in the scheme above
cannot be identified without further comments.
 To compare them, it is necessary first
to choose an equivariant isomorphism among corresponding bundles.
Then
 they coincide  up to a constant multiple.

To compare the operators $D_j$ in the above
scheme with the higher spin Dirac operators (see Def.1),
we shall  choose a certain identification of the
corresponding source and target bundles.
 We shall do it for the first operator $D_j$ in the row.

Let us characterize an algebraic operator
$Y:\Gamma(\Omega^{k+1}_c\otimes S_{\pul})\rightarrow
\Gamma(\Omega^{k}_c\otimes S_{\pul})$
by a local formula
$$
Y(\om\otimes s)=-\sum_i\iota(e_i)\om\otimes e_i\cdot s,
$$
where $\{e_i\}$ is a (local)  orthonormal basis of
$TM$ and $\iota$ denotes the contraction of a differential form
by a vector. As shown in \cite{VSe}, the map
$Y:E^{k+1,j}\rightarrow E^{k,j}, j<k<[n/2]$ is an isomorphism.

The twisted Dirac operator $D^T$ maps the space $\Gamma(\Omega^{k}_c\otimes
S_{\pul})$
to itself.
In \cite{VSe}, it was proved that we have a relation
$\nabla\circ Y+Y\circ\nabla=-D^T.$
Let us denote the projection from
$\Omega^{k}_c\otimes S$ onto $E^{k,j}$  by $\pi_{k,j}.$
 Symbols $\tilde{D}_j,\,0\leq j<[n/2]$ will denote
operators
$$
\tilde{D}_j=Y\circ D_j=
\pi_{j,j}\circ Y\circ \nabla^S|_{E^{j,j}},
$$
 mapping the space of sections of $E^{j,j}$ to itself.
Then $Y|_{E^{j,j}}=0$ implies that
$$
\tilde{D}_j=
\pi_{j,j}\circ Y\circ\nabla^S|_{E^{j,j}}=
-\pi_{j,j}\circ D^T|_{E^{j,j}}=-D_{\la_j},
$$
where $D_{\la_j}$ is the generalized higher spin Dirac operator corresponding
to the bundle $V_{\la_j},\,
\la_j=({\frac{3}{2}},\ldots,{\frac{3}{2}},{\frac{1}{2}},\ldots,
{\frac{1}{2}})$ (component ${\frac{3}{2}}$ appearing $j$
times).
More precisely, there are no signs in odd dimensions, while in
even dimension, $V_{\la_j}=V_{\la^+_j}\oplus V_{\la^-_j}.$

\section{Definition and general properties of Rarita-Schwinger operator.}
In this section, the attention is concentrated to the simplest
special case of the general definition given in above, i.e.
to the Rarita-Schwinger operator.
This operator is mentioned in \cite{RaS} and was studied in more details
 in the
paper by M.Y.Wang(\cite{Wa}) and the dissertation of U.Semmelman (\cite{US})
and also in the physical context e.g in \cite{Fra1,E,Pe1,Pe2,Pe3}.
It is also related
to the Stein-Weiss operators studied by Branson(\cite{Br}.

In the paper of Wang, it is studied
mainly on compact Einstein spin manifolds admitting nonzero
Killing spinors and in the context of deformation
of Einstein metrics on manifolds.

In the Dissertation of U.Semmelmann, the results of Wang are
extended and also some eigenvalue problems
as the relations between eigenvalues of the Dirac operator $D$,
the twisted Dirac operator $D_T$ and Rarita-Schwinger operator
$D_{\frac{3}{2}}$ are  added.

Let us first describe (a slight modification of)
 the Wang definition and his decription of the basic properties
of operators involved (\cite{Wa}).

We shall use in the following section the notation
$$
{\bf S}_{\frac{k}{2}} = V_{\la}\;\; {\rm for}\; \la =
(\frac{k}{2},\pul,...,\pul).
$$
Let ${\bf e} = \{e_1,...,e_n\}$ be an orthogonal frame on some open set
on $M$ and $\{\epsilon^1,...,\epsilon^n\}$ the dual coframe.

\begin{remark}
The classical Dirac operator on M is the operator
$$
D : \Gamma(S_{\frac{1}{2}}) \rightarrow \Gamma(S_{\frac{1}{2}})
$$
and has with respect to frame  ${\bf e}$ the form
$$
D(\sigma) = \sum_i e_i\; \nabla^S_{e_i}\sigma
$$
\end{remark}

There is a fundamental diagram, which presents the relation among
operators in discussion.

From the theory of representation of the group $Spin(n)$ we have
a decomposition of the tensor product
$$
{\bf S}_{\frac{1}{2}}\otimes \bC^{n*} = {\bf S}_{\frac{1}{2}}'\oplus
{\bf S}_{\frac{3}{2}}
$$
into irreducible
Spin(n)-modules, where ${\bf S}_{\frac{1}{2}}'$ is isomorphic to
${\bf S}_{\frac{1}{2}}$.

We shall use the following canonical identifications

$$
{\bf \iota} : {\bf S}_{\frac{1}{2}} \rightarrow
{\bf S}_{\frac{1}{2}}\otimes \bC^{n*}
$$
$$
{\bf \iota} (\sigma) :=-\frac{1}{n} \sum_i e_i\; \sigma \otimes
\epsilon^i
$$
which  is the isomorphisms onto ${\bf S}_{\frac{1}{2}}'$ and
$
{\bf S}_{\frac{3}{2}}\equiv {\rm Ker}(\mu)
$
where
$$
\mu : {\bf S}_{\frac{3}{2}}\otimes \bC^{n*}\rightarrow {\bf S}_{\frac{1}{2}}
$$
is the Clifford multiplication
$$
\mu (\sum_i \psi_i\otimes\epsilon^i) = \sum_i e_i\; \psi_i
$$
Moreover we have the projections
$
\pi_{\frac{1}{2}} , \pi_{\frac{3}{2}}
$
from
$
{\bf S}_{\frac{3}{2}}\otimes \bC^{n*}
$
onto the individual irreducible components
given by
$$
\pi_{\frac{1}{2}}(\sum_i\psi_i\otimes\epsilon^i) = -\frac{1}{n}\sum_j
e_j\;(\sum_i e_i\; \psi_i)\otimes\epsilon^j
$$
and
$$
\pi_{\frac{3}{2}}(\sum_i\psi_i\otimes\epsilon^i) = \sum_j (\psi_j
+ \frac{1}{n} e_j\; \sum_i
e_i\;\psi_i)\otimes\epsilon^j
$$
For the corresponding operators we get
\vskip 2cm
{\center{
%TexCad Options
%\grade{\off}
%\emlines{\off}
%\beziermacro{\on}
%\reduce{\on}
%\snapping{\off}
%\quality{2.00}
%\graddiff{0.01}
%\snapasp{1}
%\zoom{1.50}
\unitlength 1.00mm
\linethickness{0.4pt}
\begin{picture}(73,50)

\put(20,57){$\Gamma(S_{\frac{1}{2}})$}
\put(55,57){$\Gamma(S_{\frac{1}{2}}\otimes T_C^*)$}
\put(100,57){$\Gamma(S_{\frac{1}{2}}\otimes T_C^*)$}

\put(60,35){$\Gamma(S_{\frac{1}{2}})$}
\put(105,35){$\Gamma(S_{\frac{1}{2}})$}

\put(60,6){$\Gamma(S_{\frac{3}{2}})$}
\put(105,6){$\Gamma(S_{\frac{3}{2}})$}

\put(80,58){\vector(1,0){17.5}}
\put(75,36){\vector(1,0){23.5}}
\put(75,9){\vector(1,0){23.5}}

\put(34,58){\vector(1,0){18.5}}
%\put(77,43){\vector(1,0){17.5}}
%\put(77,28){\vector(1,0){17.5}}

\put(75,14){\vector(2,1){27}}
\put(75,31){\vector(2,-1){27}}

\put(27,53){\vector(2,-3){28}}
\put(28,54){\vector(2,-1){25}}

\put(65,45){$\|$}
\put(107,45){$\|$}

\put(65,20){$\oplus$}
\put(107,20){$\oplus$}

\put(43,62){\makebox(0,0)[cc]{$\nabla$}}
\put(89,62){\makebox(0,0)[cc]{$D_T$}}
\put(89,40){\makebox(0,0)[cc]{${\cal D}'$}}
\put(88,14){\makebox(0,0)[cc]{${\cal R}$}}

\put(40,52){\makebox(0,0)[cc]{${\cal D}$}}
\put(40,40){\makebox(0,0)[cc]{${\cal T}$}}

\put(83,31){\makebox(0,0)[cc]{${\cal T}'$}}
\put(83,21){\makebox(0,0)[cc]{${\cal T}^*$}}

\end{picture}
}}

\vskip -2cm
\vskip 3cm
\noindent with
$$
\nabla(\sigma) = \sum_i \nabla^S_{e_i}\sigma\otimes\epsilon^i
$$
and the (transported) Dirac operator
$$
{\cal D}(\sigma) := \pi_{\frac{1}{2}}(\nabla^S\sigma) =
-\frac{1}{n}\sum_i e_i\; D\sigma \otimes
\epsilon^i = {\bf \iota}(D(\sigma))
$$
the Twistor operator
$$
{\cal T}(\sigma) := \pi_{\frac{3}{2}}(\nabla^S\sigma) =
\sum_j (\nabla^S_{e_j}\sigma +\frac{1}{n} e_j \; D\sigma)
\otimes \epsilon^j
$$
and another operators which with respect to the above identifications
are
$$
{\cal D}' = -{\frac{n-2}{n}} {\bf \iota}\circ D \circ {\bf \iota}^{-1}
$$
$$
{\cal T}' = \frac{2}{n} {\cal T}\circ{\bf \iota}^{-1}
$$
and if we denote the operator
$$
\delta : \Gamma(S_{\frac{3}{2}}) \rightarrow
\Gamma(S_{\frac{1}{2}})
$$
by
$$
\delta(\sum_i \psi_i\otimes\epsilon^i) = -\sum_i
\nabla_{e_i} \psi_i
$$
then we get
$$
{\cal T}^*(\psi) = 2.{\bf \iota}(\delta(\psi))
$$

\begin{theorem}(\cite{Wa})
The operator $D_T$  with respect to the decomposition
$$
\Gamma(S_{\frac{1}{2}}\otimes T_C^*) = \Gamma(S_{\frac{1}{2}})\oplus
\Gamma(S_{\frac{3}{2}})
$$
has the following form
$$\pmatrix{
\frac{2-n}{n}{\bf \iota}\circ D\circ{\bf \iota}^{-1} &
2.{\bf\iota}\circ\delta \cr
\frac{2}{n} ({\cal T}\circ {\bf\iota}^{-1}) &
{\cal R} \cr
}$$
\end{theorem}

\section{Eigenvalues on spheres.}
The spectrum of the Dirac operator on the sphere is well known
for some time already (see \cite{Baer}).
The spectra of operators $\tilde{D}_j$
on sphere were computed in (\cite{BuSo}),
the case of a general elliptic first order operators
can be found in \cite{Br1}.

\begin{lemma}
The eigenvalues of the Dirac operator on the sphere $S_n$ with
standard metric are

$$
\mu_l=\pm\left({\frac{n}{2}}+l\right);\;l=0,1,2,\ldots.
$$
with multiplicity
$$
2^{[\frac{n}{2}]} \comb{l+n-1}{l} .
$$
\end{lemma}

The main  result of the paper \cite{BuSo} is given in the following theorem.

\begin{theorem}
  Let $D_{\la_j}=-\tilde{D}_j,\,0<j< n/2,$
 be the higher spin Dirac operators defined above,
considered on the sphere $S_n$ with the standard metric.
Then their eigenvalues are :

$$
\mu^1_l=\pm\left({\frac{n}{2}}+l\right);\;l=1,2,\ldots.
$$
with multiplicity

$$
2^{[\frac{n}{2}]}\comb{n+1}{j+1}\comb{l+n}{l-1}
\frac{(n-2j)(j+1)}{(l+j)(l+n-j)}
$$
and
$$
\mu^2_l=\pm\left[{\frac{n-2j}{n-2j+2}}\left({\frac{n}{2}}+l
\right)\right];\;l=1,2,\ldots.
$$
with multiplicity
$$
2^{[\frac{n}{2}]}\comb{n+1}{j}\comb{l+n}{l-1}
\frac{(n-2j+2)j}{(l+j-1)(l+n-j+1)}.
$$

\end{theorem}
It is proved using modification of results of Branson et others in \cite{BOO}.

The  Rarita-Schwinger operator is the operator $\tilde{D}_1$ and
its eigenvalues are:

$$
\mu^1_l=\pm\left({\frac{n}{2}}+l\right);\;l=1,2,\ldots.
$$
with multiplicity

$$
2^{[\frac{n}{2}]}\comb{n+1}{2}\comb{l+n}{l-1}
\frac{2(n-2)}{(l+1)(l+n-1)}
$$
and
$$
\mu^2_l=\pm\left[{\frac{n-2}{n}}\left({\frac{n}{2}}+l
\right)\right];\;l=1,2,\ldots.
$$
with multiplicity
$$
2^{[\frac{n}{2}]}(n+1)\comb{l+n}{l-1}
\frac{(n)}{(l)(l+n)}.
$$
The second series of eigenspaces on sphere comes out as the image
by Twistor operator of the eigenspaces of the Dirac operator on
sphere.

\section{Homogeneous solutions of Rarita-Schwinger equation
on \reo.}
The solutions of the Dirac equation on the flat space are studied
traditionally in Clifford analysis. Special attention is paid to
the solutions on $\reo := {\bf R}^{n+1}-\{0\} \;$ which
are polynomial (homogeneous) of some degree k,
namely to the functions
$$
\phi : \reo \rightarrow \ps
$$
which satisfy
\vskip 2mm
\noindent (1) $D \phi = 0$
\vskip 2mm
\noindent(2) $\phi (\lambda x) = \lambda^k \phi(x)$
for all $x\in \reo.$
\vskip 2mm
Such solutions are strongly related with the eigenfunctions of the
Dirac operator on the unit sphere $S^n \subset \reo $, the
relation follows simply from the restriction of solution to the
sphere. Such functions are also
called
spherical monogenics (of degree k).
\vskip 2mm
We make a simple generalization of the problem, instead of the Dirac
operator we shall study the Rarita-Schwinger operator, which is
define on some special spinor-valued one forms on \reo and is the
next operator in the serie of first order elliptic invariant operators
on Riemannian spin manifolds.
\vskip 2mm
We would like again to study homogenic solutions of the
Rarita-Schwinger operator on \reo and their relations with the
eigenfunctions of the corresponding operators (Dirac and
the Rarita-Schwinger ) on the sphere $ S^n$. Of course the situation is
more complicated here.
\vskip 2mm

There is a possibility to study the same problem for other
higher spin operators on spinor valued forms or higher
spin Rarita-Schwinger operators.

Let me mention only the main results, the full description
of results and another related topics
will be presented in fortcoming paper \cite{Bu1}, see also \cite{BuLN}.

We shall use on \reo fixed cartesian coordinates
$ x = (x_1,...,x_{n+1})$ and also:
$$
e_i = \frac{\partial}{\partial x_i}\;,\; x = \sum_i x_i e_i \in \bR^n\;,\;
$$
and we have
$$
\epsilon^i = dx^i\;,\; \nabla_{e_i} = \frac{\partial}{\partial x_i}.
$$
Let us take
$$
\psi = \sum_1^{n+1} \psi_i\otimes dx^i \in \ps\otimes \bR^{n*}
$$
and suppose that
$$
\psi\in \Gamma(S_{\frac{3}{2}}) (\iff \mu(\psi) = \sum_i e_i\;
\psi_i = 0.)
$$
Then the Rarita-Schwinger operator ${\cal R}$ on $\psi$ has the
form
$$
{\cal R}\psi = \sum_i (D\psi_i +\frac{1}{n+1} e_i \; (\sum_k e_k\;
D\psi_k))\otimes dx^i.
$$
\begin{remark}
The spinor valued 1-form  $\psi$  satisfying the condition
$\mu(\psi) = 0$ on M is solution of Rarita-Schwinger equation iff
there is a spinor field $\phi$ on M such that
$$
D_T \psi = {\bf \iota}(\phi).
$$
\end{remark}
Let us denote by ${\cal P}_k(1)$ the space of all polynomial k-homogeneous
solutions of Rarita-Schwinger equation on  and by ${\cal
P}_k(0)$ the space of polynomial k-homogeneous solutions of the
Dirac equation on\reo.

We would like to find a good description of the space ${\cal
P}_k(1)$,
namely its decomposition into some well defined and natural
pieces.

Let us define the map
$$
{\cal L} : \Gamma(S_{\frac{3}{2}}) \rightarrow
\Gamma(S_{\frac{1}{2}})
$$
by
$$
{\cal L}(\psi) := \sum_i x_i \psi_i
$$
Then the map is homomorphism of the corresponding bundles, namely
$$
{\cal L}(\psi + \phi) = {\cal L}(\phi) + {\cal L}(\phi),\;
{\cal L}(f.\psi) = f.{\cal L}(\psi) {\rm for} f\in C^{\infty}.
$$
\begin{lemma}
Suppose  $\psi\in {\cal P}_k(1)$, then
$$
D^3 ({\cal L}(\psi)) = 0.
$$
\end{lemma}

Let $\psi_0$ be a (k-1) homogeneous solution of the Dirac
equation on \reo, then k-homogeneous solution $\psi$ of
the Rarita-Schwinger equation satisfying
$$
D_T \psi = \sum_i e_i\; \psi_0 \otimes dx^i \hskip 3cm({\rm
eq})
$$
can be constructed in the following way.
Let
$$
\Xi : \Gamma(S_{\pul}) \rightarrow \Gamma(S_{\pul}\otimes
\Lambda^1)
$$
be the map defined by:
$$
\Xi(\psi_0) = \frac{1}{2.(n+k)} (\|x\|^2 {\cal T}(\psi_0) + \sum_j x_j.\psi_0\otimes dx^j
+ \sum_j e_j\;
(x.\psi_0)\otimes dx^j)
$$
\begin{lemma}
$\Xi(\psi_0)$ is a k-homogeneous solution of the equation
(eq).
\end{lemma}

\begin{theorem}
Let $\psi$ be a k-homogeneous solution of the Rarita-Schwinger
equation.
Then there is a unique decomposition of
$\psi$ into (k-homogeneous) pieces
$$
\psi = \psi_1 + \psi_2 + \psi_3
$$
with $\psi_j\in {\cal M}^j$
where
$$
{\cal M}^1 = \{ \psi \in \Gamma(\reo,S_{\frac{3}{2}})\;|\;  {\cal
L}(\psi) = 0 \}
$$
$$
{\cal M}^2 =  {\cal T}'(Ker D) \subset \{ \psi \in
\Gamma(\reo,S_{\frac{3}{2}})\;|\;
D{\cal L}(\psi) = 0 \}
$$
$$
{\cal M}^{3} = \Xi({\cal P}_{k-1}(0))\subset \{ \psi \in
\Gamma(\reo,S_{\frac{3}{2}})\;|\;
D{\cal L}(\psi) \neq 0\; {\rm for}\; \psi\neq 0 \}.
$$

\end{theorem}
\vskip 2mm

We can also describe the space ${\cal P}_k(1)$ of k-homogeneous solutions of
Rarita-Schwinger equation from representation point of view.
Any space ${\cal M}^j$ is an irreducible representation of the
group Spin(n) and its type is determined by  its highest weight
$\la$.
We shall speak about the space with highest weight $\la$ as about
the space of representation type $\la$.

\vskip 1mm

\noindent Recall that the Rarita-Schwinger equation for
$\psi\in \Gamma(\reo, S_{\frac{3}{2}})$ has a form:
$$
D_T \psi = \iota (\phi),
$$
with
$$
D \phi = 0.
$$
The space ${\cal P}_k(1)$  is a direct sum of three
spaces (representation types) of solutions, namely
\vskip 1mm

\vskip 3mm
\noindent${\cal M}^1$ :
\vskip 1mm
The space is characterized by the condition ${\cal L}(\psi)= 0$.
It corresponds (by restrictions of fields) to the eigenspace of induced
Rarita-Schwinger operator on unit sphere
of type $(\frac{2k+1}{2},\frac{3}{2},\pul,...,\pul)$.
\vskip 3mm
\noindent${\cal M }^2$:
\vskip 1mm
The space is characterized by the condition ${\cal L}(\psi)\neq
0\;\;{\rm for}\;\; \psi \neq 0,\; D_0{\cal L}(\psi) = 0$,
The corresponding space ${\cal P}_k^{B1}(1)$ is constructed from
the space $\{\sigma\in {\cal P}_{k+1}(0)\}$ using the Twistor operator
and is of the type $(\frac{2k+3}{2},\pul,...,\pul)$.
\vskip 2mm

Both preceding types are characterized by condition $\phi = 0$ and consists of
solutions non only of RS-equation but of the whole twisted Dirac
equation.
\vskip 2mm
\noindent${\cal M}^3$:
\vskip 1mm
The space is
characterized by the condition $D_0(\psi) \neq 0 \; {\rm for} \;
\psi\neq 0$,
and can be constructed from
${\cal P}_{k-1}(0) $
the space of (k-1)-homogenic solutions of Dirac equation,
and is of the type $(\frac{2k-1}{2},\pul,...,\pul).$

\begin{remark}
All types of solutions can be uniquely determined  by its
restriction to unit sphere and can be described by pure spherical
data (see\cite{BuLN}).
\end{remark}
\begin{remark}
The classification of polynomial solutions of Rarita-Schwinger
equation can be done
also  for polynomial solutions of the general higher spin Dirac
equation as well as for higher Rarita-Schwinger equations.
There are several papers in preparation for publication by the authors
P. Van Lancker (Gent), F. Sommen (Gent), V.Sou\v cek (Prague)
and myself.
\end{remark}

\section{Index of elliptic differential operator.}
Let $E\rightarrow M$ and $F\rightarrow M$ be complex vector
bundles over a compact m-dimensional manifold M. Denote
$\Gamma(E)$ the space of smooth sections of $E$.

Let
$$
D : \Gamma(E) \rightarrow \Gamma(F)
$$
be an elliptic differential operator between $E$ and $F$, and
$$
D^* : \Gamma(F)\rightarrow \Gamma(E)
$$
its adjoint. Both operators $D, D^*$ have finite dimensional
kernels.

The  index  of the operator $D$ is defined as
$$
{ \rm Ind} D = {\rm dim \hskip 1mm Ker}D - {\rm dim\hskip 1mm
Ker}D^*.
$$
There are formulas,  characterizing the index of an operator
in terms of topological invariants of the manifold $M$ and
the corresponding bundles.

In the formulas appear some characteristic classes of complex
vector bundle $E$ on $M$, which can be
expressed through closed differential forms representing Chern
classes $c_i(E)$  namely :

\noindent Todd genus

$$
td(E) = 1 + \frac{1}{2} c_1(E) + \frac{1}{12}(c_2(E) + c_1(E)^2)
+ ...
$$
and

\noindent  Chern character

$$
Ch(E) = {\rm dim}(E) + c_1(E) + \frac{1}{2}(c_1(E)^2 - 2 c_2(E) +
\frac{1}{6}(c_1(E)^3 - 3 c_1(E)c_2(E) + 12 c_3(E)) +...
$$
\begin{remark}
Chern classes $c_i(E)$ are represented by closed differential
forms  of degree 2i, computable from curvature form of  any
connection on $E$.

From the index theorem it follows, that index of the operator
does not depend on operator itself but only on bundles $E$ and
$F$.

In the following section, suppose that the dimension of M is
$m = 2n$  even, even if will be not especially mentioned,
because odd dimensional cases are trivial.

\end{remark}

For the index of Dirac operator and its twisted version there are
well-known index formula, containing the invariant
$\hat{A}$-genus.

\begin{theorem} (Atiyah-Singer,\cite{BGV}).
Let M be an oriented spin manifold of
dimension m = 2n. Let ${\cal S}$ be the spin bundle and let ${\cal W}$
be an arbitrary vector bundle over M, If $D_W $ is the twisted Dirac
operator on $\Gamma(M,{\cal W\otimes S})$ then
\[
ind D_W = (2\pi{\bf i})^{-n}\int_M \hat{A}(M) ch({\cal W}).
\]
\end{theorem}
We shall use normalization of forms representing characteristic
classes as is given in \cite{BGV}.

The invariant $\hat{A}(M)$ is the Hirzebruch's $\hat{A}$ genus,
given by the formula
$$
\hat{A}(M) = 1 - \frac{1}{24} p_1(M)+ (\frac{7}{5760} p_1(M)^2 -
\frac{1}{1440} p_2(M)) + ...
$$
where
$$
p_i(M) = (-1)^i c_{2i}(T_{\bf C}(M))
$$
are the Pontrjagin
forms.
\begin{remark}
Let us remark, that index  of the classical Dirac (called also
Atiyah-Singer \cite{LM} ) operator
$$
D^+_{\frac{1}{2}} : \Gamma(S^+) \rightarrow \Gamma(S^-)
$$
is
$$
Ind D^+_{\pul} := \hat{A}_{\pul}[M] = (2\pi {\bf i})^{-n}\int_M
\hat{A}(M).
$$
\end{remark}

The Chern character $Ch(E)$ of the bundle $E$ satisfies
the conditions
$$
Ch(E\oplus F) = Ch(E) + Ch(F)
$$
$$
Ch(E\otimes F) = Ch(E).Ch(F)
$$
Now we shall try to use general theory for computation of the
index of our operators $ D_j$ in discussion, namely
$$
{\rm ind} D_j^+ =  {\rm dim\; ker} D^+_j - {\rm dim\; ker} D^-_j
$$
First of all we have from the Th.1 :

\begin{corollary}
The index of the twisted Dirac operator on $M$, dim $M$ = 2n,
$$
D_T^+ : \Gamma(T^*_C(M)\otimes S^+) \rightarrow \Gamma(T^*_C(M)\otimes S^-)
$$
is given by the
formula
$$
{\rm Ind} D^+_T = \int_M Ch(T^*_C(M)) \hat{A}(M) =
$$
$$
(2\pi {\bf i})^{-n}\int_M (2n + p_1(M) + \frac{1}{12}(p_1(M)^2 - 2 p_2(M)) + .)
(1 - \frac{1}{24} p_1(M) + \frac{1}{ 5760} (7 p_1(M)^2 - 4
p_2(M)) +.).
$$

\end{corollary}

\subsection{Index of the Rarita Schwinger operator.}

Because of equality
$$
S_{\frac{1}{2}}^{\pm}\otimes T^*(M) = S_{\frac{1}{2}}^{\mp}\oplus
S_{\frac{3}{2}}^{\pm}
$$
we have
$$
Ch(S_{\frac{3}{2}}^+) - Ch(S_{\frac{3}{2}}^-) = (Ch(T^*_{\bf C}(M) +1)
(Ch(S_{\frac{1}{2}}^+) -Ch(S_{\frac{1}{2}}^-)).
$$
We can simply compute
$$
(Ch(T^*_{\bf C}(M)) + 1 = 2n+1 + p_1(M) + \frac{1}{12}(p_1(M)^2 - 2
p_2(M)) +....
$$
and finally we get
$$
{\rm Ind} D_{\frac{3}{2}} := \hat{A}_{\frac{3}{2}}[M] = (2\pi {\bf i})^{-n}\int_M
(Ch(T^*_c(M))+1)\hat{A}(M).
$$

Especially for dimension $2n = 4$ we have
$$
{\rm Ind} D_{\frac{3}{2}} =
$$
$$
(2\pi {\bf i})^{-2}\int_M (5 + p_1(M)
+...)(1-\frac{1}{24} p_1(M) + ...)
$$
$$
= (2\pi {\bf i})^{-n}\int_M (5 + \frac{19}{24}
p_1(M)) = -19 (2\pi {\bf i})^{-2}\int_M \hat{A}(M) .
$$
So we have the relation between the index of Dirac operator and
Rarita-Schwinger operator by
$$
{\rm Ind}D_{\frac{3}{2}} = -19 .{\rm Ind} D_{\frac{1}{2}}.
$$

Next nontrivial case is  dimension $m = 2n = 8$,  we have :
$$
{\rm Ind} D_{\frac{3}{2}} =
$$
$$
= (2\pi {\bf i})^{-4}\int_M (9 + p_1(M)
+\frac{1}{12} (p_1(M)^2 - 2p_2(M))+.)(1-\frac{1}{24} p_1(M)
+\frac{1}{5760} (7 p_1(M)^2 - 4 p_2(M)) + .) =
$$
$$
= (2\pi {\bf i})^{-4}\int_M (9 + \frac{19}{24}
p_1(M) + \frac{1}{5760} (543 p_1(M)^2 - 996 p_2(M)) +.).
$$
So we do not have the simple relation between the index of Dirac operator and
Rarita-Schwinger operator as above, only we have
$$
{\rm Ind}D_{\frac{3}{2}} = 249 {\rm Ind}\; D_{\pul}- \frac{21}{144}
(2\pi {\bf i})^{-4}\int_M p_1(M)^2.
$$
There is a possibility to find a manifold of dimension $*$
without harmonic spinors, but with nontrivial kernel of
Rarita-Schwinger operator.

\subsection{Index of the Dirac higher spin operators.}

For the computation of the index for operator $D_j^+$  for
$2\leq j < n $  there is an induction procedure.

We have
$$
S^{\pm} \otimes \Lambda^j T^*_C(M) \simeq V_0^{\mp}\oplus V_1^{\pm}
\oplus...\oplus V_j^{\mp} \; {\rm j \; even}
$$
$$
\hskip 2cm \simeq V_0^{\pm} \oplus V_1^{\mp}\oplus...\oplus
V_j^{\pm} \; {\rm j \;odd}.
$$
and
$$
Ch(S^{\pm}) . Ch(\Lambda^j T_C^*(M))= \sum_{k=0}^j Ch(V_k^{\mp sgn (k)})
$$
for $j$ even,
$$
Ch(S^{\pm}) . Ch(\Lambda^j T_C^*(M))= \sum_{k=0}^j Ch(V_k^{\pm sgn (k)})
$$
for $j$ odd,
with $sgn(k) = \pm $ if $k$ is even or odd.

Together we get
$$
Ch(V_j^+) - Ch(V_j^-) = (-1)^{j+1} (Ch(\Lambda^j T^*_C)
+ Ch(\Lambda^{j-1}T^*_C)). (Ch(S^+) - Ch(S^-)).
$$
\begin{theorem}
The index of the operator $D_j$ is :
$$
{\rm ind} D_j := \hat{A}_{\frac{j}{2}}[M] = (2\pi {\bf i})^{-n}. \int_M
(Ch(\Lambda^{j-1})-Ch(\Lambda^j))
\hat{A}(M).
$$
\end{theorem}
Using the theorem we can study problem of existence of
solutions of equation $\tilde{D}_j \psi = 0$.

\markright{References}

\noindent Author address: Jarol\'{\i}m Bure\v s, Mathematical Institute of
Charles University, Sokolovsk\' a 83, 186 75 Prague.


\begin{thebibliography}{99}








\bibitem[Baer]{Baer}: Baer, Ch. : {\em The Dirac operator on
space forms of positive curvature}, J.Math.Soc.Japan 48,
69-83, 1996.

\bibitem[BFGK]{BFGK} Baum, H.; Friedrich, T.; Grunewald, R.; Kath, I.: {\em
        Twistor and Killing spinors on Riemannian manifolds},
        Seminarbericht 108, Humboldt University, Berlin, 1990.

\bibitem[BGV]{BGV} Berline N., Getzler E.,Vergne M. : {\em Heat
Kernels and Dirac Operators}, Grundlehren der Math.Wiss. 298,
Springer 1992

\bibitem[BOO]{BOO} Branson T, Olafsson G., Orsted B. : {\em Spectrum
generating operators and intertwining operators for
representations induced from a maximal parabolic subgroup},
J.Funct.Analysis 135, 1996 163-205.


\bibitem[Bra]{Br} Branson T. : {\em Stein-Weiss operators and
elipticity.}, J.Funct.Anal. 151(1997) 334-383.

\bibitem[Br1]{Br1} Branson T. : {\em Spectra of self-gradients on
spheres.}
Preprint 1998.


\bibitem[BuSo]{BuSo}  Bure\v s J., Sou\v cek V.: {\em Eigenvalues of
conformally invariant operators on spheres.}, to appear in Proc.
WS 1998.

\bibitem[Bu]{Bu1} Bure\v s J. {\em Spherical monogenic 1-forms },
preprint 1998 (in preparation for publication in Proc. of Summer
School, Cetraro)

\bibitem[BuLN]{BuLN}  Bure\v s,J. : {\em The Rarita-Schwinger
operator and spherical monogenic forms} , Lectures presented at
the
University of Gent, November 1998, preprint.



\bibitem[DSS]{DSS} Delanghe, R.; Sommen, F.; Sou\v cek, V.: {\em
        Clifford Algebra and Spinor-Valued Functions,} Kluwer Ac.
        Publishers, 1992.






\bibitem[E]{E} Esposito G.:{\em Dirac operator and spectral
geometry,} Preprint, hep-th 9704016.

%\bibitem[EP]{EP} Esposito G.; Pollifrone G.: {\em Twistors and 3/2
%potentials in quantum gravity,} Twistor Newslwtter, 35-53.


\bibitem[F]{F} Fegan H.D.: {\em Conformally invariant first order
differential operators},
Quat.J.Math.Oxford , 27 (1976) 371-378.

\bibitem[Fra1]{Fra1} Frauendiener J.: {\em Another view at the spin
(3/2) equation,}
Twistor Newsletter, 37 (1994), 7-9.


\bibitem[Fra2]{Fra2} Frauendiener J.: {\em A higher spin generalization
of the Dirac equation to arbitrary curved manifolds,}
Twistor Newsletter, 37 (1994), 10-13.


\bibitem[FraS]{FraS} Frauendiener J.; Sparling G.A.: {\em On a class of
consistent linear higher spin equations on curved manifolds,}
Preprint, 1994.

\bibitem[Fri1]{Fri1} Friedrich, T.: {\em Dirac-Operatoren in der
        Riemannschen Geometrie,} Vieweg, 1997.


\bibitem[LM]{LM} Lawson H.B., Michelsom M.-L.: {\em Spin
Geometry},
Princeton Math.Series 38, Princeton Univ.Press , 1989


\bibitem[MP]{MP} Mason, L.J; Penrose, R: {\em Spin 3/2 fields and
local twistors,} Twistor Newsletter 37, (1994), 1-6.

\bibitem[NGRVN]{NGRVN} Nielsen N.K., Grisaru M.T., Roemer H.,
Van Nieuwenhuizen P. : {\em Approaches to the gravitational
Spin-$\frac{3}{2}$ axial anomaly,} Nuclear Physics B140 (1978)
477-498.


%\bibitem[N]{N} Nicolas J.-P.: {\em Spin 3/2 zero rest-mass fields in
%the Schwarzschild space-time,}
%Twistor Newsletter, 39 (1995), 6-10.

%\bibitem[NGRN]{NGRN} Nielsen N.K.; Grisaru M.T.; R"mer H.;
%Nieuwenhuizen P.: {\em Approaches to the gravitation spin-3/2
%axial anomaly,}
%Nucl.Physics, B140 (1978), 477-498.

\bibitem[Pe1]{Pe1} Penrose, R: {\em A twistor-topological approach to
the Einstein equations,} Twistor Newsletter, 38 (1994), 1-9.

\bibitem[Pe2]{Pe2} Penrose, R: {\em Twistors as spin 3/2 charges,}
Gravitation and Cosmology (A.Zichini, eds.), Plenum Press, New
York, 1991.


\bibitem[Pe3]{Pe3} Penrose, R: {\em Twistors as spin 3/2 charges
continued: SL(3,C) bundles,}
Twistor Newsletter, 33 (1991), 1-7.

\bibitem[Pe4]{Pe4} Penrose, R: {\em Twistors as charges for  3/2 in
vacuum,}
Twistor Newsletter, 32 (1991), 1-5.



\bibitem[Pe5]{Pe5} Penrose, R: {\em Concerning Space-Time points for
spin 3/2 twistor space,}
Twistor Newsletter, 39 (1995), 1-5.


\bibitem[RaS]{RaS} Rarita W.; Schwinger J.:{\em On a theory of
particles with half-integral spin,} Phys.Rev., 60 (1941), 61.


\bibitem[US]{US} Sem\-mel\-mann, U.: {\em Komp\-le\-xe
        kon\-takt\-struk\-tu\-ren und
        Kah\-ler\-sche Kil\-ling\-spi\-no\-ren,} Dissertation,
        Humbold University,
        Berlin.

\bibitem[VSe]{VSe} Severa, V.: {\em Invariant differential
        operators on Spinor-valued differential forms,}
         Dissertation, Charles University, Prague, 1998.

\bibitem[Som1]{Som1}Sommen F. : {\em Clifford analysis in two and
several vector variables}, preprint 1998.



\bibitem[So1]{So1} Sou\v cek V. : {\em Monogenic differential forms and BGG
resolution,} accepted to Proc. ISAAC Conf., Delaware, 1997.

\bibitem[SoS]{SoS}
Sommen,F; V.Sou\v cek, V: {\em Monogenic differential forms,}
 Complex Variables, Theory and Appl., 19 (1992), 81-90.

\bibitem[So2]{So2} Sou\v cek, V.: {\em Monogenic forms on
manifolds,}
        in Z. Oziewicz et. al. (Eds.), Spinors, Twistors, Clifford
        Algebras and Quantum Deformations,
         Kluwer Academic Publishers, 1993, 159-166.

\bibitem[So3]{So3}
Sou\v cek, V: {\em Conformal invariance of higher spin
equations,}
in Proc. of Symposium "Analytical and numerical methods
 in Clifford analysis,
 Seiffen, 1996.



%\bibitem[T]{T} Townsend P.K.: {\em Cosmological constant in
%supergravity,} Phys.Rev. D, 15 (1977), 2802-2805.


\bibitem[Wa]{Wa} Wang, M.: {\em Preserving Parallel Spinors under
        Metric Deformations.} Indiana Univ.Math.Jour.,
        40 (1991),




















\end{thebibliography}
\end{document}